\documentclass[11pt]{article}
\usepackage{graphics}
\usepackage{color}
\usepackage[all]{xy}

\newfont{\frak}{eufm10 scaled\magstep1}
\newfont{\sfrak}{eufm8 scaled\magstep1}
\newfont{\bbb}{msbm10 scaled\magstephalf}
\newfont{\sbbb}{msbm7 scaled\magstephalf}

\def\Ker{\hbox{Ker}}

\def\D{\Delta}
\def\C{\mbox{\bbb{C}}}
\def\R{\mbox{\bbb{R}}}
\def\Z{\mbox{\bbb{Z}}}

\def\K{\mbox{\bbb{K}}}

\def\vz{\underline{z}}

\def\n{\mbox{\frak n}}

\def\zd{(z_1,\cdots,z_d)}
\def\lorw{\longrightarrow}
\def\SC{\mbox{\sbbb{C}}}

\def\G{\Gamma}

\def\squareforqed{\hbox{\rlap{$\sqcap$}$\sqcup$}}
\def\qed{\ifmmode\else\unskip\quad\fi\squareforqed}
\def\smartqed{\def\qed{\ifmmode\squareforqed\else{\unskip\nobreak\hfil
\penalty50\hskip1em\null\nobreak\hfil\squareforqed
\parfillskip=0pt\finalhyphendemerits=0\endgraf}\fi}}

\newtheorem{thm}{Theorem}[section]
\newtheorem{prop}[thm]{Proposition}
\newtheorem{lemma}[thm]{Lemma}

\newcommand{\proof}{\mbox{\textbf{ Proof.\ \ }}}
\title{\sc Betti numbers of the geometric spaces associated to nonrational
simple convex polytopes}
\author{\sc Fiammetta Battaglia}
\date{Preliminary version, 25/04/07}
\date{}
\begin{document}
\maketitle
\begin{abstract}
\let\thefootnote\relax\footnotetext{ 
Research partially supported by MIUR (``Geometria Differenziale e Analisi
Globale'' PRIN 2007)}
We compute the Betti numbers of the geometric spaces associated to nonrational 
simple convex polytopes and find that they depend on the combinatorial 
type of the polytope exactly as in the rational case. This shows that the combinatorial features of the 
starting polytope are encoded in these generalized toric spaces as they are in their rational counterparts.
\end{abstract}

{\small 2000 \textit{Mathematics Subject Classification.} Primary: 14M25.
Secondary: 52B05, 32S99.}


\parindent0pt

\section*{Introduction}
The construction of toric varieties from rational 
convex polytopes was generalized to arbitrary  convex polytopes
in successive papers: the simple case was treated in
the joint article with Elisa Prato \cite{cx} and the nonsimple case 
in \cite{f}. The construction stems from E. Prato's paper \cite{p}, 
where the simple case is treated in the symplectic set up.
The next natural step is to explore the link
between combinatorics of convex polytopes and topological invariants of these
generalized toric varieties.
In the present article we consider the case of simple convex polytopes: 
we introduce the notion of deRham cohomology for the corresponding spaces and we compute their Betti numbers.
We find that they depend on the combinatorial type of the polytope exactly as in the rational case.
We can therefore infer that generalized toric varieties carry the same information on the combinatorics of 
the polytope as classical toric varieties, at least in the simple case.

A dimension $n$ convex polytope $\D\subset(\R^n)^*$ 
is said to be {\em simple} if there are 
exactly $n$ facets, namely codimension $1$ faces, meeting at each vertex. 
The fan dual to the polytope $\D$ is the one generated by the 
$1$--dimensional cones dual to the facets. The polytope is said to be {\em rational} if
there exists a lattice $L$ that contains  a set of generators of these $1$--dimensional dual cones.
It is known that, for each nonrational simple polytope, we can find a rational polytope of the 
same combinatorial type, namely a polytope which has the same 
$h$--vector, whilst this is not true for nonsimple polytopes (M. Perles see \cite{gru},\cite{ziegler}).
The $h$--vector $(h_0,\ldots,h_n)$ of an $n$--dimensional simple polytope 
is defined to be $$h_k=\sum_{i=0}^{k}(-1)^{k-i}{{n-i}\choose{n-k}}f_{n-i}$$
where $f_j$ is the number of faces of dimension $j$. The $h$--vector
can also be defined for nonsimple polytopes, but, in this case, it cannot be
recovered from the $f$--vector $(f_0,\ldots,f_{n})$.

Given a nonrational simple convex polytope $\D\subset(\R^n)^*$, there always exists a
{\em quasilattice} -- a $\Z$--submodule of $\R^n$ generated by a finite
set of spanning vectors of $\R^n$ -- that contains a set of generators 
of the $1$--dimensional dual cones. In \cite{cx} it was shown that, to each such a polytope, together with a 
choice of generators of the dual fan and of a quasilattice $Q\subset\R^n$
containing such generators, there corresponds a space $X$
with features that mimic those of rationally smooth toric varieties.
More precisely, $X$ is an $n$--dimensional compact complex {\em quasifold}, acted on holomorphically by the 
complex {\em quasitorus} $\C^n/Q$, with a dense open orbit; we say that $X$ is a quasitoric space.
Quasifolds and quasitori were introduced by E. Prato in \cite{p}; the notion of quasifold 
generalizes that of orbifold: quotients by discrete, not necessarily finite, groups, are allowed.
Therefore a quasifold is not necessarily Hausdorff: by compact we will mean a space such that any open covering 
admits a finite subcovering. Quasitorus is a natural generalization  of torus in this set up.  

To start with, we define the deRham cohomology of compact quasifolds and prove that, for these spaces, 
the Poincar\'e Lemma holds. Therefore, if two compact quasifolds have the same homotopy type in the $C^{\infty}$ 
sense, then they have the same deRham cohomology. Then we carry on by considering the case 
of quasitoric spaces.  
The computation of their Betti numbers is based on the following key remark: whilst 
rationally smooth toric varieties are covered by affine subsets of the kind $\C^n/\G$, where $\Gamma$  is a 
finite subgroup of the compact torus $T^n$, the quasitoric spaces corresponding to nonrational simple convex polytopes 
are covered by sets of the kind $\C^n/\G$, where $\Gamma$ is a finitely generated subgroup of $T^n$; however, 
for the computation of the Betti numbers, the only relevant fact is that $\G$ is a subgroup of the compact 
torus $T^n$. We find that, if $X$ is a quasitoric space associated to a simple convex polytope $\D$, with
$h$--vector $(h_0,\ldots,h_n)$, then the Betti numbers of $X$ are given by:
$$b_{2j+1}(X)=0\qquad\hbox{and}\qquad
b_{2j}(X)=h_j,\quad j=0,\ldots,n$$
which is exactly the relationship between Betti numbers and $h$--vector in the case of a toric variety
associated to a rational simple convex polytope.
\section{Preliminaries}
Quasifolds and related geometrical objects were introduced by E. Prato in \cite{p}, 
a recent version can be found in the joint paper with E. Prato \cite{kite}, were some of the
definitions were reformulated. 
The definition of complex quasifold is just a natural generalization of that of real quasifold, 
a version based on \cite{p} can be found in \cite{cx}.

Now we recall, from \cite{cx}, what is the quasitoric space corresponding to a simple convex polytope. 
Let $\Delta\subset (\R^n)^*$
be a simple convex polytope of dimension $n$ having $d$ facets.
We choose $d$ vectors $X_j$, that generate the $1$--dimensional dual cones, and
a quasilattice $Q$ in the $\R^n$ containing such vectors.
Once the vectors $X_j$'s and the quasilattice $Q$ have been chosen, 
the corresponding quasitoric space $X$ is the topological quotient
of an open subset, $C^d_{\D}$, of $\C^d$ by the action of a subgroup, $N_{\SC}$, of $T^d_{\SC}$, 
that can be nonclosed. Let us see in detail how  $C^d_{\D}$ and $N_{\SC}$ are defined.
We write the polytope $\D$ as
$
\D=\bigcap_{j=1}^d\{\;\mu\in(\R^n)^*\;|\;\langle\mu,X_j\rangle\geq\lambda_j\;\}.
$
Each (open) face $F$ of the polytope of dimension $p$ defines a subset 
$I_F$ of $\{1,\ldots,d\}$ such that 
$$
F=\{\,\mu\in\Delta\;|\;\langle\mu,X_j\rangle=\lambda_j\;
\hbox{ if and only if}\; j\in I_F\,\}.
$$
Since the polytope is simple the set $I_F$ contains exactly $p$ indices.
When there is no risk of ambiguity we denote $I_F$ simply by $F$.

Now let $\K$ be either $\C$ or $\C^*$.
Let $J$ be a subset of $\{1,\dots,d\}$ and let $J^c$ be its complement.
We denote by
$
\K^J=\{\zd\in\C^d\;|\;z_j\in\K\;\;\hbox{if}\;\; j\in J,\;
z_j=0\;\;\hbox{if}\;\; j\notin J\}.
$
By $T^J$ we denote the subtorus $\{(t_1,\dots,t_d)\in T^d\;|\:
t_j=1\;\hbox{if}\;j\notin J\}$. 
Consider the open subset
$$\C^d_{\D}=\cup_{F\in\D}
\C^{F}\times(\C^*)^{F^c}=\cup_{F\in\D}(\C^*)^{F^c}$$
Notice that $C^d_{\D}$ depends only on the combinatorics of the polytope.
The construction of the group $N$, due to E. Prato \cite{p},  
generalizes to non rational simple polytopes the construction by
T. Delzant \cite{delzant}.
Let $N$ be the subgroup of $T^d$ given by \linebreak 
$\Ker(\Pi:T^d\lorw \R^n/Q)$, where $\Pi$ is the mapping  
induced by the projection $\pi\,\colon\,\R^d\rightarrow\R^n$, defined by
$\pi(e_j)=X_j$, $\{e_j,\;j=1,\ldots,d\}$ being the canonical basis of $\R^d$.
The Lie algebra of $N$ is $\n=\Ker(\pi)$. 
The group $N_{\SC}$ is the complexification of $N$ and its polar decomposition
is $N_{\SC}=(\exp{i\n})N$.

The orbit space $\C^d_{\D}/N_{\SC}$, endowed with the quotient topology, is our quasitoric space $X$.
A complex quasifold structure is defined on $X$ as follows:
let $\nu$ be a vertex of the polytope $\D$ and
let $$V_{\nu}=\C^{\nu}\times(\C^*)^{\nu^c}/N_{\SC}.$$
This is an open subset of $X$.
Notice that, since $\D$ is simple, the subsets $V_{\nu}$ give an open 
covering of $X$.  
Let $\G_{\nu}$ be the finitely generated 
subgroup of $T^{\nu}$ given by $$\G_{\nu}=N\cap T^{\nu};$$ the quotient
$\C^{\nu}/\G_{\nu}$ is a  quasifold model. Consider now the the homeomorphism 
$\phi_{\nu}:\C^{\nu}/\G_{\nu}\lorw V_{\nu}$  
defined by $$\phi_{\nu}([\vz])=[\vz+\underline{1}_{\nu^c}]$$ 
with $(\underline{1}_{\nu^c})_j=1$ if $j\notin I_{\nu}$ and $(\underline{1}_{\nu^c})_j=0$ otherwise.
The triple $$\left(V_{\nu},\phi_{\nu},\C^{\nu}/\G_{\nu}\right)$$
is a {\em chart} of the space $X$; biholomorphic {\em changes of charts} are defined for each pair
of charts $V_{\nu}$ and $V_{\mu}$ (they  always have nonempty intersection).
Therefore the collection of charts indexed by the vertices of $\D$ 
is a quasifold atlas of $X$.
The quasifold $X$ is compact: this is a consequence of \cite[Thm~3.2]{cx}, stating
that there is an equivariant diffeomorphism between the complex quasitoric space $X$ and its symplectic
counterpart, constructed in \cite{p}, which is easily seen to be compact.

The general definition of quasifold goes along the above lines:
an $n$--dimensional complex quasifold is a topological space admitting 
an open covering of charts $U_{\alpha}$, such that each $U_{\alpha}$ is homeomorphic 
to the quotient of a complex manifold $\tilde{U}_{\alpha}$ modulo the  
action, by biholomorphisms, of a discrete group $\G_{\alpha}$, such that the set of points
where the action is not free has minimal real codimension $\geq2$. Moreover biholomorphic changes of charts are defined for
each pair of intersecting charts.

Recall that differential forms and the differential $d$ on $X$ are defined as follows:
a $k$ differential form $\alpha$ on $X$ is given by a collection, indexed by the vertices of 
$\D$, of differential $k$ forms
$\tilde{\alpha}_{\nu}$ on $\C^{\nu}$, invariant by the action of $\G_{\nu}$ and well behaved
under the changes of charts.
An explicit example of $2$--form is exhibited in \cite{cx}.
The differential $d$ is then defined naturally: $d\alpha$ is the $k+1$ form
given locally by $d\tilde{\alpha}_{\nu}$. We denote by $\Omega^k(X)$ the space of
$k$--differential forms on $X$.
\section{The de Rham cohomology of compact quasitoric spaces}
Let $X$ be a quasifold.
Consider the deRham complex on $X$:
$$\ldots\Omega^k(X)\stackrel{d}{\lorw}\Omega^{k+1}(X)\stackrel{d}{\lorw}
\Omega^{k+2}(X)
\stackrel{d}{\lorw}\ldots.$$
We define the {\em $k$--th deRham cohomology} of $X$ to be  $$H^k(X)=\{\hbox{closed forms in}\;\Omega^k(X)\}/
\{\hbox{exact forms in}\;\Omega^k(X)\} .$$ 
Now let $X$ and $Y$ be two quasifolds, the product
$X\times \R$ is still a quasifold, with the natural quasifold structure
induced by that of $X$ and by the differentiable structure of $\R$. 
Let $f$ and $g$ be two smooth mappings from $X$ to $Y$ (see \cite[Definition~A.33]{kite}),
$f$ and $g$ are homotopic if the usual condition is satisfied.  
We then define the homotopy type of a quasifold in the $C^{\infty}$ sense
\begin{prop} [deRham Cohomology and Homotopy Type]\label{omotopia}{\rm
Let $X$ be a quasifold and let $Y$ be of the same homotopy type as $X$,
then $X$ and $Y$ have the same deRham cohomology.}
\end{prop}
\proof
The proof of the Poincar\'e Lemma for manifolds (see for example \cite[Chap.~1]{bott-tu}), 
goes through word by word.
The homotopy operator $K$ can be defined from that on 
manifolds, since our forms are, locally, invariant forms on manifolds and 
the homotopy operator maps invariant forms to invariant forms. This gives
$H^{\bullet}(X\times \R)=H^{\bullet}(X)$. Then the invariance of the deRham cohomology under homotopy type is a 
straightforward consequence.
\qed

We now want to compute the Betti numbers of a quasitoric space $X$ associated to
the simple convex polytope $\D$ by means of the Mayer--Vietoris sequence. Before proceding with
the proof we need a preliminary lemma.
Let $M$ be a manifold endowed with the smooth  action of a group $S$. We denote the complex 
of the differential forms on $M$, that are invariant under the action of $S$, by $\Omega^{\bullet}(M)^S$. 
\begin{lemma}[Cohomology of invariant forms]\label{isomorfo}{\rm  Let $T$ be a compact connected group 
acting smoothly on a manifold $M$ and let $S$ be a subgroup of $T$, not necessarily closed. Then
$$H^k(M)\simeq H^k(\Omega^{\bullet}(M)^S)\simeq H^k(\Omega^{\bullet}(M)^T)$$}
\end{lemma}
\proof
Define the mapping $\hbox{Av}\,\colon\,\Omega^{\bullet}(M)\rightarrow\Omega^{\bullet}(M)^T$ as follows:
$$\hbox{Av}(\alpha)=\frac{1}{\hbox{Vol}(T)}\int_Tt^*(\alpha)dt.$$
Now remark that
$$\hbox{Av}\circ\hbox{Av}=\hbox{id}$$
and
$$\hbox{Av}\circ d=d\circ\hbox{Av}.$$
Moreover, since $T$ is connected, each $t\in T$ is a diffeomorphism of $M$ homotopic to the identity, therefore
$[\alpha]=[t^*(\alpha)]$ for each $t\in T$, which implies $$[\alpha]=[\hbox{Av}(\alpha)].$$
It is now easy to deduce that the natural inclusion
$\Omega^{\bullet}(M)^T\rightarrow\Omega^{\bullet}(M)$
induces an isomorphism
$$H^k(\Omega^{\bullet}(M)^T)\rightarrow H^k(M).$$
The same argument applies to the inclusion 
$\Omega^{\bullet}(M)^T\rightarrow\Omega^{\bullet}(M)^S$, yielding the isomorphism
$$H^k(\Omega^{\bullet}(M)^T)\rightarrow H^k(\Omega^{\bullet}(M)^S).$$
Therefore the mapping $\iota$ in the below diagram is an isomorphism
$$
\xymatrix{
H^k(\Omega^{\bullet}(M)^T)\ar@{->}[r]\ar@{->}@/_1.5pc/[rr]&
H^k(\Omega^{\bullet}(M)^S)\ar@{->}[r]^{\iota}& H^k(M)}
$$
\qed

In order to compute the Betti numbers of $X$, that is the dimension of the $H^k(X)$'s, we adapt to our case 
a classical argument used for the computation of the Betti numbers of rationally smooth toric varieties:
\begin{thm}[Betti numbers]{\rm Let $b_j(X)$, $k=0,\ldots,2n$, be the Betti numbers of a quasitoric space $X$ 
associated to the simple convex polytope $\D$.
Then $$b_{2k+1}=0$$ and $$b_{2k}=h_k.$$
where $h=(h_0,\cdots,h_n)$ is the $h$--vector of $\D$.}
\end{thm}
\proof
We start with some combinatorics on the simple convex  polytope $\D$.
Fix a vector $X_0\in \R^n$ such that
the hyperplanes $\langle \mu,X\rangle=c$, $c\in\R$, are transversal
to the polytope faces and do not meet more than one vertex at a time. 
Consider the height function along the direction $X_0$. Using the language of Morse theory 
we could say that the vertices of the 
polytope are the critical points of our function, we order them from the 
lowest to the highest: $\nu_1,\ldots,\nu_{f_0}$; 
to each vertex $\nu_k$ there corresponds  
the face of the polytope, $F_k$, of highest dimension among those that have
$\nu_k$ as lowest vertex. Let $r_k$ be the dimension of this face, then
$\hbox{ind}(\nu_k)=n-r_k$ is the index of the vertex $\nu_k$ as a critical point.
Remark now that $$f_0=\sum_{k=0}^{n}h_k,$$
moreover it is easy to check that the number of vertices of given index $i$ is $h_{i}$; for example:
the number of vertices of index $0$ is $h_0=f_n=1$, which corresponds to 
the only vertex of index $0$, namely the lowest vertex $\nu_{1}$;  
the number of vertices of index $1$ is $h_1=f_{n-1}-n f_{n}$, that is the total number of facets, $f_{n-1}$, minus
the $n$-facets which have $\nu_0$ as lowest vertex; the number of vertices of index $n$ 
is $h_n=\sum_{i=0}^{n}(-1)^{n-i}f_{n-i}=1$, which corresponds to 
the only vertex of index $n$, namely the highest vertex $\nu_{f_0}$.
Recall that, given a face $F$ of $\D$, we denote the subset $I_F$ of $\{1,\ldots,d\}$ simply by $F$. 
In this sense $\nu_j$ has cardinality $n$ for each $j=1,\ldots,d$, $F_k$ has cardinality $n-r_k=\hbox{ind}(\nu_k)$ and
$F_k\subset \nu_k$.
We have defined an open covering of $X$, indexed by the 
polytope vertices: we set $$V_k=V_{\nu_k}\simeq \C^{\nu_k}/\G_{\nu_k}.$$
Consider now, in $X$, the open sets 
$$W_k=\cup_{i=1}^{k}V_i,$$
each of the $W_j$'s is a quasifold (it is a union of charts) and
$$W_1\subset W_2\subset\cdots\subset W_{f_0}$$
with $W_1=V_1$ and $W_{f_0}=X$.
Recall that to each $r$--dimensional face of the polytope there corresponds an $r$--dimensional
complex quasifold, which is an orbit of the complex quasitorus $\C^n/Q$ \cite{f}.
It is easy to realize that the complement of $W_{k-1}\cap V_k$ in $V_k$  
is the union of the $\C^n/Q$--orbits associated to those faces of the polytope whose closure contain 
the vertex $\nu_k$ and that, in turn, lie in the closure
of the face $F_k$. More precisely, $V_k\setminus W_{k-1}\cap V_k$ is the quotient,
by the $\G_{\nu_k}$--action, of the $r_k$--dimensional subspace $\C^{\nu_k\cap F_K^c}$.
Therefore, as a subset of $V_k$,
\begin{equation}\label{intersezioni}
W_{k-1}\cap V_k\simeq\left((\C^{F_k}\setminus \{0\})\times\C^{\nu_k\cap F_K^c}\right)/\G_{\nu_k}.
\end{equation}
Recall that a differential form on the quasifold $W_k$ is given, on each chart $V_j$, $j=1,\ldots,k$,
by a $\G_{\nu_j}$--invariant form on $\C^{\nu_j}$.
By Lemma~\ref{isomorfo} we have that 
$$H^{\bullet}(W_{k-1}\cap V_k)\simeq
H^{\bullet}(\C^{F_k}\setminus \{0\}).$$
Now, in order to make use of the Mayer--Vietoris sequence for the computation of the Betti numbers,
we need to construct a partition of unity for the open covering $W_{k-1}\cup V_k$ of $W_k$, for
$k=2,\ldots,f_0$.
Consider a $C^{\infty}$ function $\lambda$ on $\R$ such that $0\leq\lambda\leq1$, with
$$\lambda(x)=\left\{
\begin{array}{l}
1\quad\hbox{if}\quad |x|<1\\
0\quad\hbox{if}\quad |x|>2 
\end{array}\right.$$
Let  $[\vz_1,\vz_2]\in V_{\nu_k}\simeq (\C^{F_k}\times\C^{\nu_k\cap F_k^c})/\G_{\nu_k}$ and let
$\lambda_{V_k}([\vz_1,\vz_2])=\lambda(|\vz_1|)$, this defines a smooth function on $V_k$, that
extends as a smooth function on the whole of $W_k$ (\cite[Definition~A.33]{kite}).
Consider the functions $\lambda_{V_k}$ and $\lambda_{W_{k-1}}=1-\lambda_{V_k}$ defined on
$W_k$. We remark that
the closure in $W_k$ of the set $\{x\in W_k\;|\; \lambda_{V_k}(x)\neq0\}$ is contained in $V_k$. 
Analogously the closure in $W_k$ of the set $\{x\in W_k\;|\; \lambda_{W_{k-1}}(x)\neq0\}$ is contained
in $W_{k-1}$.
Therefore $\lambda_{V_k}$ has support in $V_k$, whilst
the function $\lambda_{W_{k-1}}=1-\lambda_{V_k}$ has support in $W_{k-1}$.
The pair of functions
$$\rho_{V_k}=\frac{\lambda_{V_k}}{\lambda_{V_k}+\lambda_{W_{k-1}}}$$
$$\rho_{W_{k-1}}=\frac{\lambda_{W_{k-1}}}{\lambda_{V_k}+\lambda_{W_{k-1}}}$$
gives a partition of unity of the covering $W_{k-1}\cup V_k$ of $W_k$. This allows us to apply the argument
given in \cite[Chap.~1]{bott-tu} to constructing the Mayer--Vietoris sequence for $$W_k=W_{k-1}\cup V_k.$$ 
It is easy to prove by induction that
$$H^{2j+1}(W_{k})=0$$
for every $j$. Furthermore we find:
$$H^{2j}(W_k)=H^{2j}(W_{k-1}),\quad\quad\hbox{if}\quad j
\neq n-r_k$$
$$H^{2j}(W_k)=H^{2j}(W_{k-1})\oplus\R,\quad\quad\hbox{if}\quad j
=n-r_k.$$
This finally implies the desired result, since, as recalled above,
the number of vertices of given index $i$ is exactly $h_{i}$.
\qed 

We can conclude that the combinatorial features of the starting
polytope are encoded in the quasitoric space $X$, as it is the case 
for classical toric varieties. 

\noindent \small{\sc Dipartimento di Matematica Applicata,
Via S. Marta 3, 50139 Firenze, Italy}

\noindent \small{\sc e-mail: fiammetta.battaglia@unifi.it}

\end{document}